\documentclass[11pt]{amsart}
\usepackage{amssymb}
\usepackage[all,cmtip]{xy}

\theoremstyle{definition}

\theoremstyle{remark}

\begin{document}

\title[The Borsuk-Ulam theorem for $\mathbb{Z}_p$-actions]{A simple proof of the Borsuk-Ulam theorem for $\mathbb{Z}_p$-actions}
\author{Mahender Singh}
\address{School of Mathematics, Harish-Chandra Research Institute, Chhatnag Road, Jhunsi, Allahabad 211019, INDIA.}
\email{mahen51@gmail.com}


\begin{abstract}
In this note, we give a simple proof of the Borsuk-Ulam theorem for $\mathbb{Z}_p$-actions. We prove that, if $S^n$ and $S^m$ are equipped with free $\mathbb{Z}_p$-actions ($p$ prime) and $f: S^n \to S^m$ is a $\mathbb{Z}_p$-equivariant map, then $n \leq m$.
\end{abstract}

\maketitle

\section*{Introduction}
Let $S^n$ be the unit $n$-sphere in $\mathbb{R}^{n+1}$. There is a natural involution on $S^n$, called the antipodal involution and given by $x \mapsto -x$. The well known Borsuk-Ulam theorem states that: If there is a map $f: S^n \to S^m$ taking a pair of antipodal points to a pair of antipodal points, then $n \leq m$.

Over the years, there have been several generalizations of the theorem in many directions. We refer the reader to an interesting article \cite{Steinlein} by Steinlein, which lists 457 publications concerned with various generalizations of the Borsuk-Ulam theorem. Also the recent book \cite{Matou} by Matou\v{s}ek contains a detailed account of various generalizations and applications of the Borsuk-Ulam theorem. There are several proofs of this theorem in literature, in fact, most algebraic topology texts contains a proof.

The purpose of this note is to give a simple proof of a generalization of this theorem in the setting of group actions. Let $G$ be a group acting on a space $X$ with the action $G \times X \to X$ denoted by $(g,x) \mapsto gx$. There is associated with the group action, the orbit space $X/G$, obtained by identifying all the points in the orbit of $x$ (denoted by $\overline{x}$) for each $x\in X$. The orbit map $X \to X/G$ is given by $x \mapsto \overline{x}$.

If spaces $X$ and $Y$ carry $G$-actions, then a map $f:X \to Y$ is called $G$-equivariant if $f(gx)= g(f(x))$ for all $x \in X$ and $g \in G$. An equivariant map $f:X \to Y$ induces a map $\overline{f}: X/G \to Y/G$ given by $\overline{f}(\overline{x})=\overline{f(x)}$. Recall that a $G$-action is said to be free if $gx=x$ implies $g=e$, the identity of $G$.

In 1983, Liulevicius \cite{Liulevicius} published the following generalization of the Borsuk-Ulam theorem:

\noindent \textit{If a map $f:S^n \to S^m$ commutes with some free actions of a non-trivial compact Lie group $G$ on the spheres $S^n$ and $S^m$, then $n \leq m$.}

An alternative, but relatively simple, proof of the later theorem was also given by Dold \cite{Dold} in 1983. There are also some other generalizations of the result, see for example \cite{Biasi}. In this note, we give a simple proof of the above result for free actions of the cyclic group $\mathbb{Z}_p$ of prime order $p$ involving only elementary algebraic topology. More precisely, we prove the following theorem.\\

\noindent \textbf{Borsuk-Ulam Theorem.}
\textit{Let $S^n$ and $S^m$ be equipped with free $\mathbb{Z}_p$-actions. If there is a $\mathbb{Z}_p$-equivariant map $f: S^n \to S^m$, then $n\leq m$.}\\

Before proceeding to prove the theorem, we recall the universal coefficient theorem for singular cohomology.\\

\noindent \textbf{Universal Coefficient Theorem.} \cite[p.243]{Spanier}
\textit{There is a natural short exact sequence
$$0 \to Ext\big( H_{k-1}(X; \mathbb{Z}), \mathbb{Z}_p\big) \to H^k(X; \mathbb{Z}_p) \to
Hom \big( H_k(X; \mathbb{Z}), \mathbb{Z}_p \big) \to 0$$
for each $k \geq 0$.}

\section*{Proof of the Borsuk-Ulam Theorem}
Suppose that $n > m$. Let the $\mathbb{Z}_p$-actions on $S^n$ and $S^m$ be generated by $T$ and $S$ respectively. Note that the map $f:S^n \to S^m$ is $\mathbb{Z}_p$-equivariant if $f(T(x))= S(f(x))$ for all $x \in X$. Let $q_1: S^n \to S^n/T$ and $q_2: S^m \to S^m/S$ be the orbit maps which are also $p$-sheeted covering projections. We claim that $\overline{f}_{\#}: \pi_1(S^n/T) \to \pi_1(S^m/S)$ is zero. This will give a lift $\tilde{f}$ of $\overline{f}$, that is, the following diagram commutes
$$
\xymatrix{
S^n \ar[d]^{q_1} \ar[r]^{f} &S^m\ar[d]^{q_2}\\
S^n/T \ar[r]^{\overline{f}} \ar[ru]^{\tilde{f}} &S^m/S.}
$$

Since, $Ext\big( H_0(S^n/T; \mathbb{Z}), \mathbb{Z}_p \big)=0$, taking $k=1$ in the universal coefficent theorem, we have $H^1(S^n/T; \mathbb{Z}_p) \cong Hom \big( H_1(S^n/T; \mathbb{Z}), \mathbb{Z}_p \big)$. The same holds for $S^m/S$ also. By naturality of the universal coefficient formula, the map $\overline{f}:S^n/T \to S^m/S$ gives the following commutative diagram
$$
\xymatrix{
H^1(S^m/S; \mathbb{Z}_p) \ar[d]^{{\overline{f}}^*} \ar[r]^{\cong \hspace{7mm}} & Hom \big( H_1(S^m/S; \mathbb{Z}), \mathbb{Z}_p \big)\ar[d]^{\alpha \mapsto \alpha  {\overline{f}}_*}\\
H^1(S^n/T; \mathbb{Z}_p) \ar[r]^{\cong \hspace{7mm}}  &Hom \big( H_1(S^n/T; \mathbb{Z}), \mathbb{Z}_p \big).}
$$

For $p$ odd, both $n$ and $m$ are odd. It is known that for a free action of $\mathbb{Z}_p$ on a sphere $S^{2k-1}$, there are integers $n_1,...,n_k$ such that $S^{2k-1}/\mathbb{Z}_p$ is homotopy equivalent to the lens space $L^{2k-1}(p; n_1,...,n_k)$. Thus  both $S^n/T$ and $S^m/S$ are homotopy equivalent to lens spaces and have the following cohomology algebras \cite[p. 251]{Hatcher}
$$H^*(S^n/  T  ; \mathbb{Z}_p) \cong \mathbb{Z}_p[s,t]/ \langle s^2, t^{\frac{n+1}{2}} \rangle,$$
$$H^*(S^m/  S  ; \mathbb{Z}_p) \cong \mathbb{Z}_p[s_1, t_1]/ \langle s_1^2, t_1^{\frac{m+1}{2}} \rangle,$$

\noindent with $t= \beta(s)$ and $t_1= \beta(s_1)$, where $\beta $ is the mod-$p$ Bockstein homomorphism. Naturality of the Bockstein homomorphism gives the commutative diagram
$$
\xymatrix{
H^1(S^m/S; \mathbb{Z}_p) \ar[d]^{{\overline{f}}^*} \ar[r]^{\beta} &  H^2(S^m/S; \mathbb{Z}_p )\ar[d]^{{\overline{f}}^*}\\
H^1(S^n/T; \mathbb{Z}_p) \ar[r]^{\beta}  & H^2(S^n/T; \mathbb{Z}_p).}
$$
If ${\overline{f}}^*$ is non zero, then ${\overline{f}}^*(s_1)=s$. From the diagram we have ${\overline{f}}^*(t_1)=t$. But $0={\overline{f}}^*(t_1^{\frac{m+1}{2}})= {{\overline{f}}^*(t_1)}^{\frac{m+1}{2}}= t^{\frac{m+1}{2}}$, a contradiction as $n > m$. Hence ${\overline{f}}^*$ is zero in this case.

For $p=2$, both $S^n/T$ and $S^m/S$ have the homotopy type of real projective spaces and hence have the cohomology algebras \cite[p. 250]{Hatcher}
$$H^*(S^n/  T  ; \mathbb{Z}_2) \cong \mathbb{Z}_2[s]/ \langle s^{n+1} \rangle,$$
$$H^*(S^m/  S  ; \mathbb{Z}_2) \cong \mathbb{Z}_2[s_1]/ \langle s_1^{m+1} \rangle,$$

\noindent where $s$ and $s_1$ are homogeneous elements of degree one each.

If ${\overline{f}}^*$ is non zero, then ${\overline{f}}^*(s_1)=s$. But $0={\overline{f}}^*(s_1^{m+1})= {{\overline{f}}^*(s_1)}^{m+1}= s^{m+1}$, a contradiction as $n>m$. Hence, ${\overline{f}}^*$ must be zero and by the commutativity of the second diagram, the map $\alpha \mapsto \alpha {\overline{f}}_*$ is zero. From this we get ${\overline{f}}_*: H_1(S^n/  T  ; \mathbb{Z}) \to H_1(S^m/ S  ; \mathbb{Z})$ is zero. Now by naturality of the Hurewicz homomorphism $$h:\pi_1(S^n/T) \to  H_1(S^n/T; \mathbb{Z})$$ (which is an isomorphism in our case), we have the following commutative diagram 
$$
\xymatrix{
\pi_1(S^n/T) \ar[d]^{h}_{\cong} \ar[r]^{\overline{f}_{\#}} &\pi_1(S^m/S) \ar[d]^{h}_{\cong}\\
H_1(S^n/T; \mathbb{Z}) \ar[r]^{\overline{f}_*}  &H_1(S^m/S; \mathbb{Z}),}
$$
which shows that $\overline{f}_{\#}: \pi_1(S^n/T) \to \pi_1(S^m/S)$ is zero and hence the lift exists.

The commutativity of the first diagram shows that both $f$ and $ \tilde{f} q_1$ are lifts of $\overline{f} q_1$. Let $x_0 \in S^n$, then by definition of $q_2$, $$q_2 \big( f(x_0)\big)= q_2\big( Sf(x_0)\big)= q_2\big( S^2f(x_0)\big)=...= q_2\big( S^{p-1}f(x_0)\big),$$ that is, the fiber over $q_2 \big( f(x_0)\big)$ is the set $$\{ f(x_0), Sf(x_0),..., S^{p-1}f(x_0) \}.$$ Also $q_2 \big( \tilde{f} q_1 (x_0)\big)= \overline{f} q_1(x_0)= q_2f(x_0)$. Therefore, $\tilde{f} q_1 (x_0)=f(x_0)$ or $\tilde{f} q_1 (x_0)= S^if(x_0)$ for some $1 \leq i \leq p-1$. Note that in the later case we have $\tilde{f} q_1\big(T^i(x_0)\big)= \tilde{f}q_1(x_0)=S^if(x_0)= fT^i(x_0)$. Hence in either case, the lifts $f$ and $ \tilde{f}q_1$ agree at a point and therefore by uniqueness of lifting, we have $f = \tilde{f}q_1$. Now for any $ x \in S^n$, $q_1(x)= q_1T(x)$. But $\tilde{f}q_1(x)= \tilde{f}q_1T(x)=fT(x)=Sf(x) \neq f(x)$, a contradiction. Hence $n \leq m$. \hfill $\Box$
\bigskip

\textbf{Acknowledgement:} The author thanks the referee for comments which improved the presentation of the note.

\bibliographystyle{amsplain}

\begin{thebibliography}{10}
\bibitem {Biasi} C. Biasi, D. de Mattos, \textit{A Borsuk-Ulam Theorem for compact Lie group actions}, Bull. Braz. Math. Soc. 37 (2006), 127-137.
\bibitem {Dold} A. Dold, \textit{Simple proofs of some Borsuk-Ulam results}, (Proceedings of the Northwestern Homotopy Theory Conference, Evanston, Ill., 1982), Contemp. Math. 19, Amer. Math. Soc., Providence, R. I., 1983, 65-69.
\bibitem {Hatcher} A. Hatcher, \textit{Algebraic Topology}, Cambridge University Press, 2002.
\bibitem {Liulevicius} A. Liulevicius, \textit{Borsuk-Ulam theorems for spherical space forms}, (Proceedings of the Northwestern Homotopy Theory Conference, Evanston, Ill., 1982), Contemp. Math. 19, Amer. Math. Soc., Providence, R. I., 1983, 189-192.
\bibitem {Matou} J. Matou\v{s}ek, \textit{Using the Borsuk-Ulam theorem}, Lectures on Topological Methods in Combinatorics and Geometry, Springer, 2003.
\bibitem {Spanier} E. H. Spanier, \textit{Algebraic Topology}, Springer-Verlag, 1966.
\bibitem {Steinlein} H. Steinlein, \textit{Borsuk's antipodal theorem and its generalizations and applications: a survey}, Topological methods in nonlinear analysis, S\'{e}m. Math. Sup. 95, Presses Univ. Montr\'{e}al, Montreal, QC, 1985, 166-235.
\end{thebibliography}

\end{document}